# Benchmarking within a DEA framework: setting the closest targets and identifying peer groups with the most similar performances

*José L. Ruiz and Inmaculada Sirvent*

*Centro de Investigación Operativa. Universidad Miguel Hernández. Avd. de la Universidad, s/n 03202-Elche (Alicante), SPAIN*

**Abstract**

Data Envelopment Analysis (DEA) is widely used as a benchmarking tool for improving performance of organizations. For that purpose, DEA analyses provide information on both target setting and peer identification. However, the identification of peers is actually a by-product of DEA. DEA models seek a projection point of the unit under evaluation on the efficient frontier of the production possibility set (PPS), which is used to set targets, while peers are identified simply as the members of the so-called reference sets, which consist of the efficient units that determine the projection point as a combination of them. In practice, the selection of peers is crucial for the benchmarking, because organizations need to identify a peer group in their sector or industry that represent actual performances from which to learn. In this paper, we argue that DEA benchmarking models should incorporate into their objectives criteria for the selection of suitable benchmarks among peers, in addition to considering the setting of appropriate targets (as usual). Specifically, we develop models having two objectives: setting the closest targets and selecting the most similar reference sets. Thus, we seek to establish targets that require the less effort from organizations for their achievement in addition to identifying peer groups with the most similar performances, which are potential benchmarks to emulate and improve.

*Keywords*: Data Envelopment Analysis; Benchmarking; Peer identification; Target Setting.

## 1. Introduction

Managers of organizations use benchmarking for the evaluation of performance in comparison to best practices of others within a peer group of firms in an industry or sector. In the best-practice benchmarking process, the identification of the best firms makes it possible to set targets, taking into account considerations related to the similarity of performances or the effort needed for their achievement, in order for organizations to ultimately learn from them in planning improvements.

Data Envelopment Analysis (DEA) (Charnes et al., 1978) has been widely used as a benchmarking tool for improving performance of decision making units (DMUs). Some recent papers dealing with DEA and benchmarking include Adler et al. (2013), which uses network DEA, Zanella et al. (2013), which uses a model based on a novel specification of weight restrictions, Dai and Kuosmanen (2014), which combines DEA with clustering methods, Yang et al. (2015), which uses DEA to create a dynamic benchmarking system, Gouveia et al. (2015), which combines DEA and multi criteria decision analysis (MCDA), Daraio and Simar (2016), which deals with benchmarking and directional distances, Gomes Júnior et al. (2016), which uses non-radial efficiencies in a benchmarking analysis based on alternative targets, and Ghahraman and Prior (2016) and Lozano and Calzada-Infante (2018), which propose stepwise benchmarking approaches.

In DEA, the performance of DMUs is evaluated through setting targets on the frontier of the production possibility set (PPS) formed by the efficient units, which can be seen as a best practice frontier in the circumstance of benchmarking (see Cook et al. (2014) for discussions). Specifically, targets are the coordinates of a projection point of the unit under evaluation on to the efficient frontier, which result from a combination of efficient DMUs on a face of such frontier. Those efficient DMUs are typically grouped in the so-called "reference sets". Thus, as stated in Thanassoulis et al. (2008), *target settings and peer identification are largely the purpose of the DEA analysis* (beside the efficiency score). In particular, *where the DMU is Pareto-inefficient: identifying efficient DMUs whose operating practices it may attempt to emulate to improve its performance, and estimating target input-output levels that the DMU should in principle be capable of attaining under efficient operation* is an important part of the information most readily obtainable from a DEA analysis (see also Thanassoulis (2001)).

The selection of actual referents is particularly crucial in the benchmarking, because organizations need to identify a peer group in their sector or industry that represent actual performances from which to learn. We note that projections in DEA represent virtual DMUs that result from actual units by assuming some postulates, while efficient DMUs in the reference sets are actual DMUs (peers). As said before, DEA analyses provide information on both target setting and peer identification. However, DEA models are primarily concerned with target setting (e.g., the closest targets), while the identification of the peers (for purposes of benchmarking) can actually be seen as a by-product. That is, DEA models implement objectives that deal with the projection point that determines the targets, while the selection of the efficient DMUs that constitute the reference set is made without following any specific criterion for the benchmarking. In practice, the efficient DMUs in the reference sets are used for the benchmarking, for example by means of radar charts that allow us to compare the units under assessment against both targets and peers in order to evaluate performance and plan for improvements (see Thanassoulis et al., 2008). But it is obvious that a selection of peers based only on grounds of belonging to a DEA reference set, which simply means participating in the combination of efficient DMUs that lead to the projection point, does not necessarily ensure an appropriate selection of benchmarks, unless they are chosen following some additional desirable criterion for the benchmarking.

In this paper, we pay special attention to this aspect of the DEA analysis, namely the identification of peers for purposes of benchmarking. In the proposed approach, peer groups are identified as DEA reference sets associated with a given projection point (as usual), but they are selected following some additional desirable criteria for the benchmarking. Bearing in mind the discussion above, it can be argued that DEA benchmarking models should incorporate into the formulations objectives related to the selection of suitable benchmarks among peers, in addition to considering the setting of appropriate targets. Specifically, we seek an approach based on models that implement criteria of closeness between actual inputs and/or outputs and targets, on one hand, and similarity of performances between actual DMUs and peers, on the other.

This approach is developed within the context of non-radial DEA models. As stated in Thanassoulis et al. (2008), non-radial models are the appropriate instruments for target setting and benchmarking, because they ensure that the identified targets lie on the Pareto-efficient subset of the frontier. In particular, the models formulated here are in line with those that minimize the distance to the efficient frontier and set the closest targets. Closest targets minimize the gap between actual performances and best practices, thus showing DMUs the way for improvement with as little effort as possible. These models have made an important contribution to DEA as tool for the best-practice benchmarking of DMUs. In fact, there is already a large and growing body of research on setting closest targets and benchmarking. Portela et al. (2003) deals with similarity in DEA as closeness between actual inputs and outputs of units under evaluation and targets, and propose determining projection points as similar as possible to those units. Aparicio et al. (2007) solve theoretically the problem of minimizing the distance to the efficient frontier by using a mixed integer linear problem (MILP) that is used to set the closest targets. For the same purposes of target setting and benchmarking, the ideas in the latter paper are extended in Ramón et al. (2016) to models with weight restrictions, in Ruiz and Sirvent (2016) for developing a common benchmarking framework, in Aparicio et al. (2017) to oriented models, in Cook et al. (2017) for the benchmarking of DMUs classified in groups, in Ramón et al. (2018), which propose a sequential approach for the benchmarking, in Ruiz and Sirvent (2018) to models that incorporate information on goals, and in Cook et al. (2019) for a benchmarking approach within the context of pay-for-performance incentive plans[1]. In the present paper, the models developed have two objectives: minimizing the gap between actual inputs and/or outputs and targets and minimizing the dissimilarities between actual performances and those of the members of the reference sets. Thus, these models will allow us to evaluate performance in terms of targets that require the less effort for their achievement while at the same time identifying peers with the most similar performances to emulate and improve.

---

[1] We note that the idea of closeness to the efficient frontier has also been investigated for purposes of developing efficiency measures that satisfy some desirable properties. See Ando (2017), Aparicio and Pastor (2014), Fukuyama et al. (2014a, 2014b) and Fukuyama et al. (2016).

We develop non-oriented and oriented DEA benchmarking models, taking into account in addition the specification of returns to scale (variable or constant). The two objectives are implemented in most cases through the closeness between vectors of inputs and/or outputs, except for the selection of the reference set in the constant returns to scale (CRS) case, wherein the similarity of performances between actual DMUs and peers is assessed in terms of the deviations between the corresponding mixes of inputs and outputs. In all of the models, a parameter is included in the formulations that allows us to adjust the importance attached to each of the two objectives. Through the specification of such parameter, a series of targets and reference sets can be generated. These series provide managers with more flexibility in target setting and benchmark selection, both for evaluating performance and for future planning. This is illustrated through an example in which teaching performance of public Spanish universities is evaluated.

The paper is organized as follows: In section 2, we develop a bi-objective DEA non-oriented benchmarking model that allows to setting the closest targets and selecting the closest reference set. Section 3 addresses the case of oriented models (either to inputs or to outputs). Section 4 addresses the case of having a technology in which constant returns to scale is assumed. In section 5 we examine an example for purposes of illustration. Last section concludes.

## 2. A model for setting the closest targets and selecting the closest reference set

Throughout this section, we consider that we have $n$ DMUs whose performances are to be evaluated in terms of $m$ inputs and $s$ outputs. These are denoted by $(X_j, Y_j)$, $j=1,...,n$, where $X_j = (x_{1j},...,x_{mj})' > 0_m$, $j=1,...,n$, and $Y_j = (y_{1j},...,y_{sj})' > 0_s$, $j=1,...,n$. For the benchmarking, we assume a variable returns to scale (VRS) technology (Banker et al., 1984). Thus, the production possibility set (PPS), $T = \{(X,Y)/X \text{ can produce } Y\}$, can be characterized as

$$T = \left\{(X,Y) \Big/ X \geq \sum_{j=1}^{n} \lambda_j X_j, Y \leq \sum_{j=1}^{n} \lambda_j Y_j, \sum_{j=1}^{n} \lambda_j = 1, \lambda_j \geq 0\right\}.$$

For a given DMU$_0$, the DEA models that set the closest targets find a projection point $(\hat{X}_0, \hat{Y}_0)$ by minimizing the distance from that unit to the points on the efficient frontier of T, $\partial(T)$, dominating $(X_0, Y_0)$. The problem below provides a general formulation of such models

$$\text{Min} \quad d\{(X_0, Y_0), (\hat{X}_0, \hat{Y}_0)\}$$

s.t.:

$$(\hat{X}_0, \hat{Y}_0) \in \partial(T)$$

$$(\hat{X}_0, \hat{Y}_0) \text{ dominates } (X_0, Y_0) \quad (1)$$

Typically, the distance between two points of T is calculated by using the following weighted $L_1$-norm: $d\{(X_0, Y_0), (\hat{X}_0, \hat{Y}_0)\} = \|(X_0, Y_0) - (\hat{X}_0, \hat{Y}_0)\|_1^\omega = \sum_{i=1}^{m} \frac{(x_{i0} - \hat{x}_{i0})}{x_{i0}} + \sum_{r=1}^{s} \frac{(\hat{y}_{r0} - y_{r0})}{y_{r0}}$ (note that we do not use absolute values because the projection point dominates DMU$_0$). If deviations between actual inputs and outputs and the coordinates of the projection are expressed by using the classical slacks, then the model below provides an operational formulation of (1) (see Aparicio et al., 2007)

$$\text{Min} \quad \sum_{i=1}^{m} \frac{s_{i0}^-}{x_{i0}} + \sum_{r=1}^{s} \frac{s_{r0}^+}{y_{r0}}$$

s.t.:

$$\sum_{j \in E} \lambda_j x_{ij} = x_{i0} - s_{i0}^- \quad i = 1, ..., m$$

$$\sum_{j \in E} \lambda_j y_{rj} = y_{r0} + s_{r0}^+ \quad r = 1, ..., s$$

$$\sum_{j \in E} \lambda_j = 1$$

$$-\sum_{i=1}^{m} v_i x_{ij} + \sum_{r=1}^{s} u_r y_{rj} + u_0 + \delta_j = 0 \quad j \in E$$

$$v_i \geq 1 \quad i = 1, ..., m$$

$$u_r \geq 1 \quad r = 1, ..., s$$

$$\lambda_j \delta_j = 0 \quad j \in E$$

$$\delta_j, \lambda_j \geq 0 \quad j \in E$$

$$s_{i0}^-, s_{r0}^+ \geq 0 \quad \forall i, r$$

$$u_0 \text{ free} \quad (2)$$

where E is the set of extreme efficient DMUs in T.

Remark. In Aparicio et al. (2007) the constraints $\lambda_j \delta_j = 0$, $j \in E$, are equivalently expressed by using the classical big M and binary variables. Eventually, in practice, we use Special Ordered Sets (SOS) (Beale and Tomlin, 1970) for solving the resulting problem. SOS Type 1 is a set of variables where at most one variable may be nonzero. Therefore, if we remove these constraints from the formulation and define instead a SOS Type 1 for each pair of variables $\{\lambda_j, \delta_j\}$, $j \in E$, then it is ensured that $\lambda_j$ and $\delta_j$ cannot be simultaneously positive for DMU$_j$'s, $j \in E$. As a result, the DMU$_j$'s in the reference set of DMU$_0$ belong all to the same supporting hyperplane of the T at the projection point of that unit. SOS variables have already used for solving models like (2) in Ruiz and Sirvent (2016, 2018), Aparicio et al. (2017), Cook et al. (2017) and Cook et al. (2019).

Model (2) finds the closest dominating projection point to DMU$_0$ on the efficient frontier of T. The coordinates of such point allow us to set targets in terms of which the performance of DMU$_0$ can be evaluated. Specifically, these targets can be expressed by using the optimal solutions of (2) as

$$\hat{x}_{i0}^* = x_{i0} - s_{i0}^{-*} = \sum_{j \in E} \lambda_j^* x_{ij}, \quad i = 1,...,m$$
$$\hat{y}_{r0}^* = y_{r0} - s_{r0}^{+*} = \sum_{j \in E} \lambda_j^* y_{rj}, \quad r = 1,...,s$$
(3)

Model (2) allows us also to identify a reference set for DMU$_0$, $RS_0 = \{DMU_j / \lambda_j^* > 0\}$ (as usual in DEA analyses). As said in the introduction, the DMUs in $RS_0$ (actual units) are typically used for benchmarking purposes. Note, however, that $RS_0$ is identified by model (2) as a by-product of setting the closest targets to DMU$_0$, because no specific desirable criterion for the benchmarking is followed in its selection. $RS_0$ simply consists of the efficient DMUs on a face of $\partial(T)$ that determine the projection of DMU$_0$, $\left(\hat{X}_0^*, \hat{Y}_0^*\right)$, as a combination of them. As a result, in practice it often happens that peers with quite dissimilar performances to that of DMU$_0$ are selected as benchmarks for its evaluation, simply because they are members of the combination of efficient DMUs that determines the closest targets.

For these reasons, in this paper we argue that DEA benchmarking models should consider in their objectives not only the setting of appropriate targets but also the selection of suitable benchmarks among peers. Specifically, the models we develop here consider the following two objectives:

Objective 1: The targets set should be as close as possible to the actual inputs and outputs of the unit under evaluation.

Objective 2: The selected reference sets should consist of peers representing actual performances as similar as possible to that of the unit under evaluation.

Objective 1 seeks to set targets that require the least possible effort from the $DMU_0$ for their achievement. Objective 2 seeks to identify a peer group representing potential benchmarks that $DMU_0$ can emulate as easily as possible for improving performance. Model (1) only considers the first objective in its formulation. To address the second objective, $RS_0$ can be selected by minimizing the distance from $DMU_0$ to all of its members. In order to do so, we follow here a minmax approach, that is, we minimize the distance from $DMU_0$ to the unit in the reference set with most dissimilar performance. Actually, using a minmax approach for minimizing globally the distances from $DMU_0$ to all the members of $RS_0$ is equivalent to minimizing the Hausdorff distance from a point to a finite set (in a metric space). In this context, the Hausdorff distance, which we denote by $d_H\{(X_0, Y_0), RS_0\}$, is defined as

$$d_H\{(X_0, Y_0), RS_0\} =$$
$$= \max_{j \in RS_0} d\{(X_0, Y_0), (X_j, Y_j)\} = \max_{j \in RS_0} \|(X_0, Y_0) - (X_j, Y_j)\|_1^\omega = \max_{j \in RS_0} \left( \sum_{i=1}^m \frac{|x_{i0} - x_{ij}|}{x_{i0}} + \sum_{r=1}^s \frac{|y_{r0} - y_{rj}|}{y_{r0}} \right) \quad (4)$$

Ideally, one would like to minimize simultaneously both distances $d\{(X_0, Y_0), (\hat{X}_0, \hat{Y}_0)\}$ and $d_H\{(X_0, Y_0), RS_0\}$ for target setting and peer identification. However, it is obvious that the targets that minimize $d\{(X_0, Y_0), (\hat{X}_0, \hat{Y}_0)\}$ will not be necessarily the coordinates of a projection point determined by the closest peers to $DMU_0$ (that is, those that minimize $d_H\{(X_0, Y_0), RS_0\}$), and the other way around too. As a compromise between these two objectives, we propose to minimize a convex combination of them

$$\alpha d\{(X_0, Y_0), (\hat{X}_0, \hat{Y}_0)\} + (1-\alpha) d_H\{(X_0, Y_0), RS_0\}, \quad (5)$$

where $0 \leq \alpha \leq 1$. This will provide us with non-dominated solutions[2]. Through the specification of $\alpha$, we may adjust the importance that is attached to each of the two objectives. The specification $\alpha = 1$ leads to the closest targets to $DMU_0$ on $\partial(T)$, that is, those provided by model (2). However, as said before, the selection of suitable benchmarks among peers in that case is not considered as an objective in the model. As $\alpha$ decreases, the selection of closer peers is considered together with the setting of closer targets. Obviously, selecting peers with more similar performance could entail setting targets that are more demanding, although that is not always necessarily the case in real applications, as will be shown later in the empirical illustration. In practice, it may be useful to provide a series of targets and peer groups that are generated by specifying values for $\alpha$ in a grid between 0 and 1 (this approach has already been used in Stewart (2010) for benchmarking models that incorporate management goals). These series may offer different views on organizations' performance and provide managers with different alternatives that they can consider both in evaluating performance and in planning improvements. This all will be illustrated in section 5.

Bearing in mind the above, a model that sets the closest targets that are sought for $DMU_0$ while at the same time selecting the closest reference set associated with such targets (for a given $\alpha$) can be formulated as follows

$$\begin{aligned} \text{Min} \quad & \alpha d\left\{(X_0, Y_0), (\hat{X}_0, \hat{Y}_0)\right\} + (1-\alpha) d_H\left\{(X_0, Y_0), RS_0\right\} \\ \text{s.t.:} \quad & \\ & (\hat{X}_0, \hat{Y}_0) \in \partial(T) \\ & (\hat{X}_0, \hat{Y}_0) \text{ dominates } (X_0, Y_0) \end{aligned} \qquad (6)$$

In terms of the classical slacks, the following model provides an operational formulation of (6)

---

[2] To ensure non-dominated solutions in the cases α=0 and 1, a second stage would have to be carried out, in which the objective that is not considered in (5) is minimized subject to (5) equals the optimal value obtained in the first stage.

$$\text{Min} \quad \alpha\left(\sum_{i=1}^{m}\frac{s_{i0}^{-}}{x_{i0}}+\sum_{r=1}^{s}\frac{s_{r0}^{+}}{y_{r0}}\right)+(1-\alpha)z_{0}$$

s.t.:

$$\sum_{j\in E}\lambda_{j}x_{ij}=x_{i0}-s_{i0}^{-} \qquad i=1,...,m$$

$$\sum_{j\in E}\lambda_{j}y_{rj}=y_{r0}+s_{r0}^{+} \qquad r=1,...,s$$

$$\sum_{j\in E}\lambda_{j}=1$$

$$-\sum_{i=1}^{m}v_{i}x_{ij}+\sum_{r=1}^{s}u_{r}y_{rj}+u_{0}+\delta_{j}=0 \qquad j\in E$$

$$v_{i}\geq 1 \qquad i=1,...,m$$

$$u_{r}\geq 1 \qquad r=1,...,s$$

$$\lambda_{j}\delta_{j}=0 \qquad j\in E$$

$$\lambda_{j}\leq I_{j} \qquad j\in E$$

$$d\{(X_{0},Y_{0}),(X_{j},Y_{j})\}I_{j}\leq z_{0} \qquad j\in E$$

$$I_{j}\in\{0,1\} \qquad j\in E$$

$$\delta_{j},\lambda_{j}\geq 0 \qquad j\in E$$

$$s_{i0}^{-},s_{r0}^{+}\geq 0 \qquad \forall i,r \qquad (7)$$

$$z_{0}\geq 0$$

$$u_{0} \text{ free}$$

where the scalars $d\{(X_{0},Y_{0}),(X_{j},Y_{j})\}=\sum_{i=1}^{m}\frac{|x_{i0}-x_{ij}|}{x_{i0}}+\sum_{r=1}^{s}\frac{|y_{r0}-y_{rj}|}{y_{r0}}$, $j\in E$, are to be calculated prior to solving the model.

Model (7) minimizes both $d\{(X_{0},Y_{0}),(\hat{X}_{0},\hat{Y}_{0})\}$ and $d_{H}\{(X_{0},Y_{0}),RS_{0}\}$ by using a minsum approach through the parameter $\alpha$. Like in (2), the set of constraints of model (7) restricts the search of projections to points on $\partial(T)$ that dominate DMU$_0$. Note that because of the constraints $\lambda_{j}\delta_{j}=0$, $j\in E$, if $\lambda_{j}>0$ then $\delta_{j}=0$, $j\in E$. This means that if DMU$_j$ in E belongs to a reference set of DMU$_0$ then it is necessarily on the hyperplane $-v'X+u'Y+u_{0}=0$. Thus, all the DMUs in a reference set are on a same face of $\partial(T)$, because they all are on the same supporting hyperplane of T, whose coefficients are non-zero. Therefore, solving (7) leads to a projection point of DMU$_0$, $(\hat{X}_{0}^{*},\hat{Y}_{0}^{*})$, onto a face of $\partial(T)$ spanned by the DMUs in the reference set, $RS_{0}=\{DMU_{j}/\lambda_{j}^{*}>0\}$. However, in (7), $RS_{0}$ is identified by minimizing in addition the distance from DMU$_0$ to that set (in terms of the Hausdorff distance). By virtue of the constraints $\lambda_j \leq I_j$, $j \in E$, if $DMU_j$ belongs to $RS_0$ then $I_j = 1$. Note that model (7) minimizes $z_0$, which is the maximum of the distances from $DMU_0$ to all the $DMU_j$'s in E with $I_j = 1$, and consequently, it minimizes the distance from $DMU_0$ to the reference sets associated with all of the solutions of the model.

We highlight the fact model (7) provides, as usual in the standard DEA analysis, a bundle of input-output targets and the corresponding reference set, but now they all are determined following criteria of closeness to the $DMU_0$ under evaluation. Targets are the coordinates of the projection point $(\hat{X}_0^*, \hat{Y}_0^*)$, and can be obtained as in (3), in this case by using the optimal solutions of (7). This projection point results from a combination of the DMUs in the reference set selected $RS_0 = \{DMU_j / \lambda_j^* > 0\}$, which are identified by minimizing the differences between the inputs and outputs of $DMU_0$ and those of the most dissimilar peer in that set of efficient DMUs, while at the same time minimizing the differences between $DMU_0$ and the projection point. As said before, model (7) can be solved for different specifications of $\alpha$, thus generating a series of targets and peers that may provide managers with different alternatives for the benchmarking analysis.

## 3. The case of oriented models

The models developed in the previous section consider inefficiencies both in inputs and in outputs for target setting and peer identification. However, in many of the applications that are carried out in practice, performance is evaluated by using oriented models. That is, empirical applications are often intended to evaluate the potential for expanding outputs with the actual level of resources or the potential for saving resources with the actual level of outputs. In this section, we develop an oriented version of model (7).[3] The developments are made in the output oriented case (the extension to input oriented models is straightforward). Thus, the model we propose is aimed at setting the closest output targets for the level of inputs of $DMU_0$, while at the same time identifying peers having the most similar performances to that of $DMU_0$.

Setting targets in this situation means that $DMU_0$ should be projected onto one of the points of the PPS having inputs not greater than $X_0$ and outputs that cannot be improved with the actual level of inputs, $X_0$. That is, onto one point $(\hat{X}_0, \hat{Y}_0)$ of the PPS with $\hat{X}_0 \leq X_0$ whose output bundle $\hat{Y}_0$ ($\hat{Y}_0 \geq Y_0$) belongs to $\partial(P(X_0))$, which is the efficient frontier of $P(X_0) = \{Y/(X_0,Y) \in T\}$. For an appropriate target setting, $\hat{Y}_0$ should be set as the closest point in $\partial(P(X_0))$ (objective 1). In addition, since we seek benchmarks of output performance for a given level of inputs, the objective of performance similarity between $DMU_0$ and peers (objective 2) is implemented in oriented models through the closeness of both inputs and outputs (like in the non-oriented model (7)). Therefore, a general formulation of the benchmarking oriented model that is wanted is the one below

$$\text{Min} \quad \frac{\alpha}{s} d^O(Y_0, \hat{Y}_0) + \frac{1-\alpha}{m+s} d_H\{(X_0, Y_0), RS_0\}$$
$$\text{s.t.:} \quad \hat{Y}_0 \in \partial(P(X_0))$$
$$\hat{Y}_0 \geq Y_0$$

(8)

where $\alpha$, $0 \leq \alpha \leq 1$, is used for the same purposes as in (6) and $d^O(Y_0, \hat{Y}_0) = \|Y_0 - \hat{Y}_0\|_1^{\omega_O} = \sum_{r=1}^{s} \frac{\hat{y}_{r0} - y_{r0}}{y_{r0}}$.

---

[3] See Aparicio et al. (2017) for an approach based on bilevel linear programming for determining the least distance to the efficient frontier and setting the closest targets that uses an oriented version of the Russell measure.

In the case of oriented models, the distances $d^O$ and $d_H$ are divided by s and m+s, respectively (thus representing an average of deviations), in order to balance the magnitude of the two quantities that are aggregated in the objective function of the model. Note that $d^O$ only accounts for deviations in outputs while $d_H$ does it for both inputs and outputs.

In order to find an operational formulation of (8), we note that $\hat{Y}_0 \in \partial(P(X_0))$ if, and only if, the optimal value of the following problem equals zero: $\text{Max} \left\{ \sum_{r=1}^{s} s_{r0}^+ \middle| \sum_{j \in E} \lambda_j x_{ij} + t_{i0} = x_{i0}, i = 1,...,m, \right.$

$\left. \sum_{j \in E} \lambda_j y_{rj} - s_{r0}^+ = \hat{y}_{r0}, r=1,...,s, \sum_{j \in E} \lambda_j = 1, \lambda_j \geq 0, j \in E, t_{i0} \geq 0, i=1,...m, s_{r0}^+ \geq 0, r=1,...,s \right\}$. And the optimal value of this problem is zero if, and only if, there exists a feasible solution with $s_{r0}^+ = 0$, $r=1,...,s$, satisfying the optimality conditions. That is, if, and only if, there exist $t_{i0} \geq 0, i=1,...,m$, $v \geq 0_m$, $u \geq 1_s$, $\delta_j, \lambda_j \geq 0$, $j \in E$, $u_0 \in \mathbb{R}$, such that $\sum_{j \in E} \lambda_j x_{ij} + t_{i0} = x_{i0}$, $i=1,...,m$, $\sum_{j \in E} \lambda_j y_{rj} = \hat{y}_{r0}$, $r=1,...,s$, $\sum_{j \in E} \lambda_j = 1$, $-v'X_j + u'Y_j + u_0 + \delta_j = 0$, $j \in E$, $\lambda_j \delta_j = 0, j \in E$, and $v_i t_{i0} = 0, i=1,....,m$.

Therefore, the model below provides an operational formulation of (8)

$$\text{Min} \quad \frac{\alpha}{s} \sum_{r=1}^{s} \frac{s_{r0}^+}{y_{r0}} + \frac{1-\alpha}{m+s} z_0$$

s.t.:

$$\sum_{j \in E} \lambda_j x_{ij} = x_{i0} - t_{i0} \qquad i = 1,...,m$$

$$\sum_{j \in E} \lambda_j y_{rj} = y_{r0} + s_{r0}^+ \qquad r = 1,...,s$$

$$\sum_{j \in E} \lambda_j = 1$$

$$-\sum_{i=1}^{m} v_i x_{ij} + \sum_{r=1}^{s} u_r y_{rj} + u_0 + \delta_j = 0 \qquad j \in E$$

$$u_r \geq 1 \qquad r = 1,...,s$$

$$\lambda_j \delta_j = 0 \qquad j \in E$$

$$v_i t_{i0} = 0 \qquad i = 1,...,m$$

$$\lambda_j \leq I_j \qquad j \in E$$

$$d\{(X_0, Y_0), (X_j, Y_j)\} I_j \leq z_0 \qquad j \in E$$

$$I_j \in \{0,1\} \qquad j \in E$$

$$\delta_j, \lambda_j \geq 0 \qquad j \in E$$

$$s_{r0}^+ \geq 0 \qquad r = 1,...,s \qquad (9)$$

$$v_i, t_{i0} \geq 0 \qquad i = 1,...,m$$

$$z_0 \geq 0$$

$$u_0 \text{ free}$$

Targets for the outputs in terms of the optimal solutions of model (9) can be obtained as $\hat{y}_{r0}^* = y_{r0} - s_{r0}^{+*} = \sum_{j \in E} \lambda_j^* y_{rj}, r = 1,...,s$, the reference set for DMU$_0$ being $RS_0 = \{DMU_j / \lambda_j^* > 0\}$, as usual.

## 4. The constant returns to scale (CRS) case

The specification of returns to scale for the technology has some implications when the objective of performance similarity for peer identification (objective 2) is addressed in the formulation of DEA benchmarking models. In the VRS case, the two objectives established in this paper for the benchmarking models can be implemented through the closeness between the inputs and/or outputs of DMU$_0$ and both targets and the inputs and/or outputs of the peers in the reference set. However, if constant returns to scale (CRS) is assumed, the similarity of performances between DMU$_0$ and the actual efficient DMUs in the reference set cannot be measured in terms of closeness between the corresponding vectors of inputs and outputs, because in the CRS case the members of the reference set may be efficient DMUs operating at scales that are very different from that of DMU$_0$. Instead, similarity can be measured in terms of that of the corresponding mixes of inputs and outputs. Specifically, we should minimize the deviations between the mixes of inputs and outputs of DMU$_0$ and those of the DMUs in the reference set, which are measured through the sines of the corresponding angles (see Coelli (1998) and Cherchye and Van Puyenbroeck (2001) for models that identify points on the efficient frontier by minimizing mix deviations).

An alternative formulation to (7) for the CRS case is the model below

$$\text{Min} \quad \frac{\alpha}{m+s}\left(\sum_{i=1}^{m}\frac{s_{i0}^{-}}{x_{i0}}+\sum_{r=1}^{s}\frac{s_{r0}^{+}}{y_{r0}}\right)+\frac{(1-\alpha)}{2}z_0$$

s.t.:

$$\sum_{j\in E}\lambda_j x_{ij} = x_{i0} - s_{i0}^{-} \qquad i=1,...,m$$

$$\sum_{j\in E}\lambda_j y_{rj} = y_{r0} + s_{r0}^{+} \qquad r=1,...,s$$

$$-\sum_{i=1}^{m}v_i x_{ij} + \sum_{r=1}^{s}u_r y_{rj} + \delta_j = 0 \qquad j\in E$$

$$v_i \geq 1 \qquad i=1,...,m$$

$$u_r \geq 1 \qquad r=1,...,s$$

$$\lambda_j \delta_j = 0 \qquad j\in E$$

$$\lambda_j (1-I_j) = 0 \qquad j\in E$$

$$\left[m^I(X_0,X_j)+m^O(Y_0,Y_j)\right]I_j \leq z_0 \qquad j\in E$$

$$I_j \in \{0,1\} \qquad j\in E$$

$$\delta_j, \lambda_j \geq 0 \qquad j\in E$$

$$s_{i0}^{-}, s_{r0}^{+} \geq 0 \qquad \forall i,r \qquad (10)$$

$$z_0 \geq 0$$

where the scalar $m^I(X_0,X_j)$ is the sine of the angle between the input mixes of DMU$_0$ and DMU$_j$, $j\in E$. That is, $m^I(X_0,X_j) = \sin(X_0,X_j)$, where $\sin^2(X_0,X_j) = 1 - \cos^2(X_0,X_j)$,

$$\cos(X_0,X_j) = \frac{\sum_{i=1}^{m}x_{i0}x_{ij}}{\sqrt{\sum_{i=1}^{m}x_{i0}^2}\sqrt{\sum_{i=1}^{m}x_{ij}^2}}, \; j\in E.$$ Analogously, the sine of the angle between the output mixes,

$m^O(Y_0, Y_j)$, $j \in E$, can be defined. Note that these scalars have to be calculated prior to solving the model.

Model (10) works in a similar manner as (7). The main difference is in that objective 2 is now expressed in terms of similarities of input and output mixes between $DMU_0$ and the $DMU_j$'s in a reference set. In (7) we use the $L_1$-distance (weighted), which aggregates the deviations between vectors of inputs and outputs, while in (10) we also aggregate the deviations between the mixes of inputs and outputs as measured through the sines of the corresponding angles. Again, the two components of the objective function are divided by a constant in order to balance their magnitude. Specifically, the one associated with objective 1 is divided by m+s while that of objective 2 is divided by 2. The constraints $\lambda_j(1-I_j) = 0$ ensure that, if $\lambda_j > 0$, then $I_j = 1$ (these non-linear constraints can be handled by using SOS variables as explained in the Remark on model (2)). Therefore, model (10) seeks to minimize the maximum deviation between the mixes of inputs and outputs of $DMU_0$ and those of the $DMU_j$'s in the corresponding reference set. Finally, note that the convexity constraint and the variable $u_0$ have been removed as the result of assuming a CRS technology.

## 5. Illustrative example

For purposes of illustration only, in this section we apply the proposed approach to the evaluation of teaching performance of public Spanish universities. The existing regulatory framework in Spain emphasizes the importance of undergraduate education as a key issue in the performance of the universities, in line with the higher education policies and practices in Europe. This regulation also raises the need to design mechanisms for the evaluation of the universities' teaching performance, independently from their performance in other areas, such as those of research or knowledge transfer. Specifically, after the reforms that were motivated by the convergence with European Space of Higher Education (ESHE), the new degrees are required to have a quality assurance system available, which must include procedures for the evaluation and improvement of the quality of teaching and of the academic staff.

We carry out an evaluation of teaching performance from a perspective of benchmarking and target setting. The practice of benchmarking is growing among universities, which see comparisons with their peers as an opportunity to analyze their own strengths and weaknesses and to establish directions for improving performance. In this study, we pay special attention to the selection of suitable benchmarks among universities through the approach that has been proposed here, which allows us in addition to set the closest targets to the universities under evaluation.

The public Spanish universities may be seen as a set of homogeneous DMUs that undertake similar activities and produce comparable results regarding teaching performance, so that a common set of outputs can be defined for their analysis. Specifically, for the selection of the outputs, we only take into consideration variables that represent aspects of performance explicitly mentioned as requirements in the section "Expected Results" of the Royal Order RD 1393/2007, such as graduation, drop out and progress of students. As for the inputs, we consider academic staff and expenditures, as two variables that account for human and physical capital, and total enrollments, as a proxy of the size of the universities, which has a significant impact on the performance of European universities (see Daraio et al. (2015)). We adopt the academic year 2014-15 as the reference year (this approach has been followed in similar models used for performance evaluation of universities; see, for example, Agasisti and Dal Bianco (2009)). The variables considered are defined as explained below:

OUTPUTS

- GRADUATION (GRAD): Total number of students that complete the programme of studies within the planned time.
- RETENTION (RET): Total number of students enrolled for the first time in the academic year 2012-13 that keep enrolled at the university.
- PROGRESS (PROG): Total number of passed credits[4].

INPUTS

- TOTAL ENROLLMENT (ENR): Total number of enrollments.
- ACADEMIC STAFF (ASTF): Full-time equivalent academic staff.
- EXPENDITURE (EXP): This input exactly accounts for expenditure on goods and services after the item corresponding to works carried out by other (external) companies has been removed. EXP thus reflects the budgetary effort made by the universities in the delivery of their activities.

Data for these variables have been taken from the corresponding report by the Conference of Rectors of the Spanish Universities (CRUE). The sample consists of 38 (out of 48) public Spanish universities. Table 1 reports a descriptive summary.

Table 1. descriptive summary

---

[4] Credit is the unit of measurement of the academic load of the subject of a programme.

An initial DEA analysis reveals 18 universities as technically efficient: University of Almería (UAL), University of Granada (UGR), University of Huelva (UHU), University of Málaga (UMA), University of Sevilla (USE), University of La Laguna (ULL), University of Castilla La Mancha (UCLM), University Autónoma of Barcelona (UAB), University of Barcelona (UBA), University of Girona (UDG), University Pompeu Fabra (UPF), University of Extremadura (UEX), University of La Rioja (URI), University of Burgos (UBU), University Autónoma of Madrid (UAM), University Carlos III de Madrid (UC3M), University Rey Juan Carlos (URJC) and University Pública de Navarra (UPN). The reference sets used in the evaluation of the different inefficient universities consist of benchmarks selected from among those efficient universities. As has been discussed throughout the paper, the selection of benchmarks among peers can be seen as a by-product of DEA models, in particular of those that minimize the distance to the efficient frontier (like (2)), which seek to set the closest targets. Thus, these models identify the reference set that determines the projection point that establishes the targets that require less effort for improvement. However, in practice, there may be other reference sets having peers whose performances are more similar to that of the unit under evaluation, thus being more suitable benchmarks from which to learn. For example, as we shall see later, model (2) selects URI and UAM as referents in the evaluation of URV, in spite of having very dissimilar performances. This suggests that it may be worth searching for other peers that are more appropriate benchmarks. We can see in the distance matrix (Table 2), which records the distances between some inefficient universities and the efficient ones, that UAL, UDG and UPF (among others) are closer to URV than URI and UAM, thus becoming more appropriate benchmarks (potentially) as they have more similar performances to that of URV.

Table 2. distance matrix

For the evaluation of performance, we use the output-oriented model (9). This model generates a series of targets and peers by considering simultaneously the objectives of closeness between universities and both targets and peers (the specification of the parameter α allows us to attach different degrees of importance to each of the two objectives). Tables 3 and 4 report the results regarding peer identification and target setting (respectively) provided by model (9) for University Rovira i Virgili (URV), University of Zaragoza (UZA), University of Alicante (UA) and University of Valladolid (UVA), which are considered as representative cases. For each of these four universities and for a grid of α values (from 1 to 0.1 by 0.1)[5], Table 3 reports the reference sets (for each member of a reference set, its distance to the university under evaluation is recorded between brackets) and Table 4 reports the targets. The last column of each table records the distance components that are aggregated in the objective function of (9), that is, those regarding the distances from the university under evaluation both to its projection on the efficient frontier ($d^O/3$) and to its reference set ($d_H/6$). The rows α=1 actually record the results provided by model (2), which sets the closest targets on the efficient frontier. They are therefore included in the tables as the reference specification. As α decreases, model (9) attaches more importance to the selection of closer benchmarks among peers. The changes in the distance $d^O/3$ provide us with an insight into the extra effort that is needed for the achievement of the new targets, while the changes in $d_H/6$ allow us to assess the degree to which the new reference sets are more similar.

---

[5] Results exclude the specification α=0, because targets in that case are established without considering the closeness to actual inputs and outputs.

In the evaluation of URV with model (2) (α=1), the following reference set is found $R^1_{URV} = \{UDG, UPF, URI, UAM\}$. However, when α is lower than 1, model (9) identifies a new reference set $R^\alpha_{URV} = \{UDG, UPF, UC3M\}$, (α<1). Therefore, solving the model suggests removing URI and UAM from $R^1_{URV}$ (because, as said before, they are quite different from URV) and considering UC3M whose performance is more similar. Note that $d_H/6$ reduces from 0.976 to 0.346 when the reference set is changed. It is also important to highlight that the output targets associated with α=1 and those of α's lower than 1 are quite similar (see in Table 4 that $d^O/3$ only raises from 0.107 to 0.161 when changing the reference set). This means that model (9) allows us to set targets that require a similar effort as that needed with the closest targets for their achievement, albeit those targets are associated with a reference set consisting of peers whose performances are more similar to that of URV, thus being more appropriate benchmarks. The same happens with UZA, because the reference set provided by (2) includes some universities that are very dissimilar from UZA (UBA and UPN) while the one provided by (9) suggests other peers that are more like this university (UCLM and UAB). Now, $d_H/6$ reduces from 0.728 to 0.212, so, again, there would be no substantial effect on target setting ($d^O/3$ only raises from 0.051 to 0.087); perhaps, a little extra effort in GRAD should be made, but in exchange for setting a less demanding target in RET.

In the case of UA, the model when α<1 suggests removing UAL and UGR and considering UAM and URJC, which are closer. The same happens with UVA. However, in the latter case, model (9) identifies a third reference set, $R^\alpha_{UVA} = \{UCLM, UAM, UC3M\}$, when $\alpha \leq 0.8$. The values of $d_H/6$ show that new reference sets including peers with more similar performances can be selected for the benchmarking of UVA, and this flexibility entails in turn more freedom in setting targets.

Nevertheless, we should note that selecting in particular $R_{UVA}^{\alpha}$, $\alpha \leq 0.8$, for the benchmarking, would be suggesting to make an extra effort mainly in GRAD and PROG (with respect to that determined by closest targets), regardless of having a less demanding target for RET.

Table 3. Model (9) benchmarking results

Table 4. Model (9) target setting results

## 6. Conclusions

Benchmarking within a DEA framework is carried out through the setting of targets and the identification of a peer group as the reference set associated with the projection point on the efficient frontier that determines the targets. However, peer identification is actually a by-product of DEA models, whose objectives are particularly concerned with the setting of targets (for example, with finding the closest targets). The approach proposed in this paper seeks to develop DEA benchmarking models having as an objective not only the setting of appropriate targets but also the identification of peers following desirable criteria for benchmark selection. Specifically, the models formulated have two objectives: setting the closest targets and selecting the closest reference sets. The results of the empirical application show that, in practice, these models often identify reference sets consisting of peers with more similar performances to the unit under evaluation than those provided by the models that are only concerned with setting the closest targets. Furthermore, it has been observed that, in many cases, selecting peers with more similar performances does not necessarily entail setting targets that require an extra effort for their achievement. In any event, the models proposed, through the specification of parameter α, may generate a series of targets and peers that offer different alternatives for managers to consider in evaluating performance and in future planning. Potential lines of future research would include, on one hand, the extension of the models proposed to deal with units that are classified in groups whose members experience similar circumstances, following the ideas in Cook and Zhu (2007) and Cook et al. (2017) and, on the other, their extension to incorporate management goals, following, for example, the approaches in Stewart (2010) and Cook et al. (2019).

**Acknowledgments**

This research has been supported through Grant MTM2016-76530-R (AEI/FEDER, UE).

|  | ENR | ASTF | EXP | GRAD | RET | PROG |
| --- | --- | --- | --- | --- | --- | --- |
| Mean | 20362.0 | 1639.9 | 15513793.9 | 1845.2 | 3801.4 | 795103.2 |
| Standard Dev. | 12105.0 | 949.4 | 9828600.3 | 1057.1 | 2197.4 | 466860.9 |
| Minimum | 3735 | 371.9 | 2514366.7 | 399 | 755 | 140645.5 |
| Maximum | 55662 | 3855.1 | 46471378.3 | 4804 | 9442 | 1983413 |

Table 1. Descriptive summary

| Ineff./Eff | UAL | UGR | UHU | UMA | USE | ULL | UCLM | UAB | UBA | UDG | UPF | UEX | URI | UBU | UAM | UC3M | URJC | UPN |
|---|---|---|---|---|---|---|---|---|---|---|---|---|---|---|---|---|---|---|
| URV | 1.309 | 13.581 | 2.006 | 7.704 | 14.991 | 2.678 | 4.761 | 8.600 | 16.048 | 0.892 | 1.222 | 3.671 | 4.226 | 2.422 | 5.856 | 2.074 | 6.615 | 2.117 |
| UZA | 3.668 | 2.952 | 4.103 | 1.269 | 3.603 | 2.413 | 1.038 | 1.173 | 4.368 | 3.533 | 3.081 | 2.207 | 5.183 | 4.336 | 0.623 | 2.168 | 1.615 | 4.201 |
| UA | 3.167 | 4.935 | 3.664 | 1.862 | 5.736 | 1.641 | 0.404 | 2.597 | 6.865 | 2.976 | 2.752 | 1.360 | 4.998 | 3.955 | 1.114 | 1.336 | 1.909 | 3.768 |
| UVA | 3.116 | 5.131 | 3.623 | 2.084 | 5.951 | 1.561 | 0.638 | 2.700 | 7.026 | 2.927 | 2.643 | 1.272 | 4.984 | 3.929 | 1.123 | 1.233 | 2.070 | 3.741 |

Table 2. Distance matrix

|  |  | Reference set |  |  |  | Distance ($d_H/6$) |
|---|---|---|---|---|---|---|
| URV | α |  |  |  |  |  |
|  | 1 | UDG (0.892) | UPF (1.222) | URI (4.226) | UAM (5.856) | 0.976 |
|  | ≤0.9 | UDG (0.892) | UPF (1.222) | UC3M (2.074) |  | 0.346 |
| UZA | α |  |  |  |  |  |
|  | 1 | UGR (2.952) | UBA (4.368) | UAM (0.623) | UPN (4.201) | 0.728 |
|  | ≤0.9 | UMA (1.269) | UCLM (1.038) | UAB (1.173) |  | 0.212 |
| UA | α |  |  |  |  |  |
|  | 1 | UCLM (0.404) | UC3M (1.336) | UAL (3.167) | UGR (4.935) | 0.822 |
|  | ≤0.9 | UCLM (0.404) | UC3M (1.336) | UAM (1.114) | URJC (1.909) | 0.318 |
| UVA | α |  |  |  |  |  |
|  | 1 | UCLM (0.638) | UC3M (1.233) | UAL (3.116) | UGR (5.131) | 0.855 |
|  | 0.9 | UCLM (0.638) | UC3M (1.233) | UAM (1.123) | URJC (2.070) | 0.345 |
|  | ≤0.8 | UCLM (0.638) | UC3M (1.233) | UAM (1.123) |  | 0.206 |

Table 3. Benchmarking with model (9)

|  |  | Targets | | | | | | Distance |
|---|---|---|---|---|---|---|---|---|
|  |  | Inputs | | | Outputs | | | |
|  |  | ENR | ASTF | EXP | GRAD | RET | PROG |  |
| URV | α/actual | **11587** | **1071.75** | **13402500.7** | **1468** | **2079** | **518174** | $d^O/3$ |
|  | 1 | 11587 | 1071.75 | 13402500.7 | 1511.1 | 2477.5 | 570132.9 | 0.107 |
|  | ≤0.9 | 11587 | 905.0 | 13402500.7 | 1661.5 | 2460.9 | 539671.4 | 0.119 |
| UZA | α/actual | **27054** | **2834.5** | **22531331** | **2512** | **5179** | **1115003** |  |
|  | 1 | 27054 | 2216.9 | 22531331 | 2621.4 | 5746.4 | 1115003 | 0.051 |
|  | ≤0.9 | 27054 | 2149.0 | 22531331 | 3040.4 | 5324.6 | 1141579.1 | 0.087 |
| UA | α/actual | **24322** | **2036.25** | **15204775** | **2378** | **4234** | **948324** |  |
|  | 1 | 24322 | 1767.6 | 15204775 | 2378 | 4885.4 | 950313.4 | 0.052 |
|  | ≤0.9 | 24322 | 1555.0 | 15204775 | 2378 | 4822.2 | 1014266.4 | 0.069 |
| UVA | α/actual | **23249** | **2117** | **15670022.3** | **2315** | **4158** | **883332** |  |
|  | 1 | 23249 | 1670.7 | 15670022.3 | 2315 | 4791.8 | 901236.6 | 0.058 |
|  | 0.9 | 23249 | 1376.9 | 15670022.3 | 2315 | 4720.0 | 977228.3 | 0.080 |
|  | ≤0.8 | 23249 | 1880.6 | 15670022.3 | 2490.5 | 4592.4 | 991503.7 | 0.101 |

Table 4. Target setting with model (9)